\documentclass[11pt]{article}

\usepackage{amsfonts}
\usepackage{enumerate}

\usepackage{amssymb}

\makeatletter
\renewcommand{\@biblabel}[1]{#1.\hfill}
\makeatother
\usepackage{amssymb}
\begin{document}

\centerline{\bf On the Skitovich--Darmois theorem for the group of
$p$-adic numbers}

\bigskip

\centerline{\bf Gennadiy  Feldman}

\bigskip

\centerline{G.M. Feldman,  B.Verkin Institute for Low
Temperature  Physics}

\centerline{and Engineering of the National Academy of
Sciences of Ukraine,}

\centerline{
Kharkov, Ukraine feldman@ilt.kharkov.ua}

\bigskip

\begin{abstract}

Let $\Omega_p$ be the group of  $p$-adic numbers,  $ \xi_1$ and
$\xi_2$ be independent random variables with values in $\Omega_p$
and   distributions   $\mu_1$ and $\mu_2$. Let $\alpha_j, \beta_j$
be topological automorphisms of $\Omega_p$. Assuming that the linear
forms $L_1=\alpha_1\xi_1 + \alpha_2\xi_2$ and $L_2=\beta_1\xi_1 +
\beta_2\xi_2$ are independent, we describe possible distributions
$\mu_1$ and $\mu_2$ depending on the automorphisms $\alpha_j,
\beta_j$. This theorem is an analogue for the group $\Omega_p$ of
the well-known Skitovich--Darmois theorem, where a Gaussian
distribution on the real line is characterized by the independence
of two linear forms.

\end{abstract}


\bigskip


\section{Introduction}
\label{}

The classical characterization theorems of mathematical statistics
were extended to different algebraic structures such as locally
compact Abelian groups, Lie groups, quantum groups, symmetric
spaces (see e.g. \cite{Fe4}-- \cite{Fe18},
\cite{Fe20bb}--\cite{GraLo}, and also \cite{Fe6}, where one can
find necessary references). In particular,    much attention has
been devoted to  the study of the Skitovich--Darmois theorem,
where a Gaussian distribution is characterized by the independence
of two linear forms, for some classes of locally compact Abelian
groups, and the Heyde theorem, where a Gaussian distribution is
characterized by the symmetry of  the conditional distribution of
one linear form given another. In these cases coefficients of
linear forms are topological automorphisms of a group. The article
is devoted to the Skitovich--Darmois theorem for  the group of
$p$-adic numbers $\Omega_p$. To the best of our knowledge the
characterization problems for the group $\Omega_p$ have not been
studied earlier.

We recall that according to the classical Skitovich--Darmois
theorem, if $\xi_j$, $j = 1, 2, \dots, n,$ $ n \ge 2,$ are
independent random variables,  $\alpha_j, \ \beta_j $ are nonzero
constants, and  the linear forms $L_1=\alpha_1\xi_1 + \cdots
 +\alpha_n\xi_n$ and $L_2=\beta_1\xi_1 + \cdots + \beta_n\xi_n$
are independent, then all random variables $\xi_j$ are Gaussian.
This theorem was generalized by Ghurye and Olkin to the case when
$\xi_j$ are independent vectors in the space
 $\mathbb{R}^m$ and the coefficients $\alpha_j,\beta_j$ are non
 singular matrices. They proved that the independence of $L_1$ and
 $L_2$ implies that
all random vectors $\xi_j$ are Gaussian  (\cite[Ch. 3]{KLR}).

 Let $X$ be a
second countable locally compact Abelian group, ${\rm Aut}(X)$ be
the group of topological automorphisms of $X$, $\xi_j,$ $j = 1, 2,
\dots, n,$ $ n \ge 2,$ be independent random variables with values
in $X$ and   distributions $\mu_j$. Consider the linear forms
$L_1=\alpha_1\xi_1 + \cdots + \alpha_n\xi_n$ and $L_2=\beta_1\xi_1 +
\cdots + \beta_n\xi_n,$ where $\alpha_j, \beta_j \in {\rm Aut}(X)$.
In the earlier papers the main attention was paid to the following
problem: For which groups $X$ the independence of $L_1$ and $L_2$
implies that all $\mu_j$ are either Gaussian distributions, or
belong to a class of distributions which we can consider as a
natural analogue of the class of Gaussian distributions. This
problem was studied for different classes of locally compact Abelian
groups (\cite[Ch. V]{Fe6}).  It turned out that in contrast to the
classical situation, the cases of $n=2$ and an arbitrary $n$ are
essentially different. For $n=2$ this problem was solved for the
class of finite Abelian groups in \cite{Fe4}, for the class of
compact totally disconnected Abelian groups in \cite{FeGra1}, and
for the class of discrete Abelian groups in \cite{FeGra2}. We  also
note that  group analogues of the Skitovich--Darmois theorem for
$n=2$ are closely connected with the positive definite functions of
product type introduced by Schmidt (see  \cite{Sch},  \cite{Fe7}).

In the  article we continue these investigations. On the one hand,
we prove that the Skitovich--Darmois theorem, generally speaking,
fails for the group of $p$-adic numbers $\Omega_p$. On the other
hand, we give the complete descriptions of all automorphisms
$\alpha_j, \beta_j \in {\rm Aut}(\Omega_p)$  such that the
independence of the linear forms $L_1=\alpha_1\xi_1 +
\alpha_2\xi_2$ and $L_2=\beta_1\xi_1 + \beta_2\xi_2$ implies that
$\mu_1$ and $\mu_2$ are idempotent distributions, i.e. shifts of
the Haar distributions of compact subgroups of $\Omega_p$. We note
that since $\Omega_p$ is a totally disconnected group, the
Gaussian distributions on $\Omega_p$ are degenerated (\cite[Ch.
4]{P}).

\section{Definitions and notation}
\label{}

We will use some results of the duality theory for the locally
compact Abelian groups (see  \cite{HeRo}). Before we formulate the
main theorem we recall some definitions and agree on notation. For
an arbitrary locally compact Abelian
 group $X$  let $Y=X^\ast$ be its character group,
and  $(x,y)$  be the value of a character $y \in Y$ at an element $x
\in X$. If $K$ is a closed subgroup of $X$, we denote by $A(Y, K) =
\{y \in Y: (x, y) = 1$ for all $x \in K \}$ its annihilator. If
$\delta : X \mapsto X$ is a continuous homomorphism,
  then the adjoint homomorphism $\widetilde \delta : Y \mapsto Y$
is defined by the formula $(x, \widetilde \delta y) = (\delta x, y)$
for all $x \in X, \ y \in Y$. We note that $\delta \in {\rm Aut}(X)$
if and only if $\widetilde\delta \in {\rm Aut}(Y)$. Denote by $I$
the identity automorphism of a group.

Let  ${M^1}(X)$ be the convolution semigroup of probability
distributions on $X.$ For a distribution $\mu \in {M^1}(X)$ denote
by $$\widehat \mu(y) = \int_X (x, y) d \mu(x)$$ its characteristic
function (Fourier transform), and by $\sigma(\mu)$ the support of
$\mu$. For $\mu \in {M^1}(X)$, we define the distribution $\bar \mu
\in M^1(X)$ by the formula $\bar \mu(E) = \mu(-E)$ for any Borel set
$E \subset X$. Observe that $\widehat {\bar \mu}(y) =
\overline{\widehat \mu(y)}$. Let $K$ be a compact subgroup of  $X$.
Denote by $m_K$  the Haar distribution on $K$. We note that the
characteristic function of
 $m_K$ is of the form
\begin{equation}\label{2}
\widehat m_K(y) = \cases{ 1, \quad  \quad y \in A(Y,K), \cr\cr 0,
\quad  \quad y \notin A(Y,K).  \cr}
\end{equation}
Denote by $I(X)$ the set of all idempotent distributions on $X$,
i.e. the set of shifts of the Haar distributions $m_K$ of the
compact subgroups $K$ of $X$. Let $x\in X$. Denote by $E_x$ the
degenerate distribution concentrated at the point $x$.

\section{The main theorem}
\label{}

Let $p$ be a prime number. We need some properties of the group of
$p$-adic numbers $\Omega_p$ (see  \cite[\S 10]{HeRo}). As a set
$\Omega_p$ coincides with the set of sequences of integers of the
form
 $x=(\dots,x_{-n}, x_{-n+1},\dots,
x_{-1}, x_0, x_1,\dots,x_n, x_{n+1},\dots),$ where $x_n \in\{0,
1,\dots, p-1\}$, such that  $x_n=0$ for $n < n_0$, where the number
$n_0$ depends on $x$. We correspond to each element $x \in \Omega_p$ the
series
 $\sum\limits_{k=-\infty}^{\infty} x_k p^k.$ Addition and  multiplication of
the series are defined in a natural  way and they define
the operations of addition and  multiplication in  $\Omega_p$. With
respect to these operations $\Omega_p$ is a field. Denote by
$\Lambda_k$ a subgroup of $\Omega_p$ consisting of  $x \in \Omega_p$
such that $x_n=0$ for $n < k$. The subgroup $\Lambda_0$ is called
the group of $p$-adic integers and is denoted by $\Delta_p$. We note
that $\Lambda_k=p^k\Delta_p$. The family of the subgroups
 $\{\Lambda_k\}_{k=-\infty}^{\infty}$ forms an open basis at zero of
 the group $\Omega_p$ and defines a topology on $\Omega_p$.
With respect to this topology the group $\Omega_p$ is locally
compact, non-compact, and totally disconnected. We note that the
group
 $\Omega_p$ is represented as a union $\Omega_p= \bigcup\limits_{l=-\infty}^{\infty}
p^l \Delta_p$. The character group $\Omega_p^{*}$ of the group
$\Omega_p$ is topologically isomorphic to $\Omega_p$, and the value
of a character $y \in \Omega_p^{*}$ at an element $x \in \Omega_p$
is defined by the formula
 \begin{equation}\label{1}
(x, y)= \exp\Big[2\pi i\Big(\sum_{n=-\infty}^\infty
x_n\Big(\sum_{s=n}^\infty y_{-s}p^{-s+n-1}\Big)\Big)\Big],
\end{equation}
where for given $x$ and $y$ the sums in (\ref{1}) actually are
finite. Each automorphism $\alpha \in {\rm Aut}(\Omega_p)$ is of the
form $\alpha g=x_{\alpha} g,$ $g \in \Omega_p,$ where $x_\alpha \in
\Omega_p,$ $x_\alpha \ne 0$. For $\alpha \in {\rm Aut}(\Omega_p)$
 we  identify the automorphism
 $\alpha \in {\rm
Aut}(\Omega_p)$ with the corresponding element
 $x_\alpha
\in \Omega_p$, i.e. when we write $\alpha g$, we  suppose that
$\alpha \in \Omega_p$. We note that $\widetilde\alpha=\alpha$.
Denote by $\Delta_p^0$ the subset of $\Omega_p$  consisting of all
invertible in $\Delta_p$ elements, $\Delta_p^0=\{x \in \Omega_p: x_n
=0$ for $n < 0, \ x_0 \ne 0 \}$. We note that each element  $g \in
\Omega_p$ is represented in the form $g = p^k c,$ where $k$ is an
integer, and
 $c \in \Delta_p^0$. Hence, multiplication on $c$ is a
topological automorphism of the group $\Delta_p$.

Denote by ${\mathbb Z}(p^\infty)$ the set of rational numbers of
the form $\{{k / p^n} : k=0, 1, \dots,p^n-1, \ n=0,1,\dots\}$. If
we define the operation in ${\mathbb Z}(p^\infty)$ as addition
modulo 1, then ${\mathbb Z}(p^\infty)$ is transformed into an
Abelian group which we consider in the discrete topology.
Obviously, this group is topologically isomorphic to the
multiplicative group of all $p^n$th roots of unity, where $n$ goes
through the set of nonnegative integers, considering in the discrete
topology. For a fixed $n$ denote by ${\mathbb Z}(p^n)$ a subgroup
of ${\mathbb Z}(p^\infty)$ consisting of all elements of the form
${\{{k / p^n} : k=0, 1, \dots,p^n-1\}}$. Note that the group
 ${\mathbb Z}(p^n)$ is topologically isomorphic to
the multiplicative group of all $p^n$th roots of unity,
considering in the discrete topology. Observe that the groups
${\mathbb Z}(p^\infty)$ and $\Delta_p$ are the character groups of
one another.

Now we will  prove the main result of the paper. We will do this
for the linear forms
 $L_1=\xi_1 + \xi_2$ and $L_2=\xi_1 + \alpha\xi_2$,
where $\alpha \in {\rm Aut}(\Omega_p)$, and then will show how the
general case is reduced to this one.

\bigskip

{\bf Theorem 1}  {\it Let $X = \Omega_p$, $\alpha \in {\rm Aut}(X)$,
$\alpha = p^k c$, $c \in \Delta_p^0.$ Then the following statements
hold.

{$1$.} Assume that either $k=0$ or $|k|=1$. Let $\xi_1$ and $\xi_2$
be independent random variables with values in $X$ and distributions
$\mu_1$ and $\mu_2$. Assume that the linear forms  $L_1 = \xi_1 +
\xi_2$ and $L_2=\xi_1 + \alpha \xi_2$  are independent. Then

${1(i)}$ If $k=0$, then $\mu_1, \mu_2 \in I(X)$; moreover if $c=(0, 0, \dots, 0, 1, c_1, \dots)$,
 then $\mu_1$ and $\mu_2$ are degenerate distributions;

${1(ii)}$ If $|k|=1$, then either $\mu_1 \in I(X)$ or $\mu_2 \in I(X)$.

{$2$.} If $|k| \ge 2,$ then there exist independent random variables
$\xi_1$ and $\xi_2$  with values in $X$ and distributions $\mu_1$
and $\mu_2$ such that  the linear forms
 $L_1 = \xi_1 + \xi_2$ and  $L_2=\xi_1 +  \alpha \xi_2$
are independent whereas $\mu_1, \mu_2 \notin I(X)$.}

\bigskip

To prove Theorem 1 we need some lemmas. Let $\xi$ be a random
variable with  values in a second countable locally compact
Abelian group $X$ and
 distribution  $\mu$.  Taking into account that the characteristic function of the
distribution  $\mu$ is the expectation  ${\bf E}[(\xi, y)]$, exactly
as in the classical case, we may   prove the following statement.

\bigskip

{\textsc{\bf Lemma 1}} {\it Let $X$ be a second countable locally
compact Abelian group. Let $\xi_1$ and $\xi_2$ be independent random
variables with  values in $X$ and
  distributions $\mu_1$ and $\mu_2$. Then the
independence of the linear forms $L_1=\xi_1 +  \xi_2$ and $L_2=
\xi_1 + \alpha\xi_2$, where $\alpha \in {\rm Aut}(X),$ is equivalent
to the fact that the characteristic functions $\widehat\mu_1(y)$ and
$\widehat\mu_2(y)$ satisfy   the equation}
\begin{equation}\label{3}
\widehat\mu_1(u + v)\widehat\mu_2(u + \widetilde\alpha
v)=\widehat\mu_1(u)\widehat\mu_2(u)
\widehat\mu_1(v)\widehat\mu_2(\widetilde\alpha v), \quad u, v \in Y.
\end{equation}

\bigskip

{\bf Lemma 2}  {\it Let $X = \Omega_p$ and $\alpha \in {\rm
Aut}(X),$ $\alpha = p^k c,$ $c \in \Delta_p^0.$ Let
  $\xi_1$ and $\xi_2$ be independent random variables with  values in $X$ and
  distributions $\mu_1$ and $\mu_2$ such that
$\mu_j(y) \ge 0,$ $j = 1, 2.$ Assume that the linear forms
 $L_1
= \xi_1 + \xi_2$ and $L_2=\xi_1 + \alpha \xi_2$ are independent.
Then there exists a subgroup $B= p^l\Delta_p$ in $Y$ such that
$\widehat\mu_j(y)=1$ for $y \in B, \ j=1, 2$.}

\bigskip

{\it Proof}  We use the fact that the family of the subgroups $\{p^l
\Delta_p\}_{l=-\infty}^\infty$ forms an open basis at zero of
 the group $Y$. Since $\widehat\mu_1(0)=\widehat\mu_2(0)=1$, we can choose
$m$ in such a way that  $\widehat\mu_j(y) >0$ for $y \in L=p^m
\Delta_p,$ $j=1, 2$. Put $M=L$  if $k \ge 0$, and $M=p^{-k}L$  if
$k<0$. Then $M$ is a subgroup of $L$ and
 $\alpha(M) \subset L$.
Put $\psi_j(y)= - \log\widehat\mu_j(y),$ $y \in L,$ $j=1, 2$.

By Lemma 1 the   characteristic functions $\widehat\mu_j(y)$
satisfy equation (\ref{3}). Taking into account that
$\widetilde\alpha=\alpha$, we get from (\ref{3}) that the
functions $\psi_j(y)$ satisfy the equation
\begin{equation}\label{4}
\psi_1(u+v)+\psi_2(u+\alpha v)=\psi_1(u)+\psi_2(u)+
\psi_1(v)+\psi_2(\alpha v), \quad u \in L, \ v \in M.
\end{equation}
Integrating equation (\ref{4}) over the group $L$ with respect to
the Haar distribution $dm_L(u)$ and using the fact that the Haar
distribution $m_L$ is $L$-invariant, we obtain
$$
\psi_1(v)+\psi_2(\alpha v)=0, \quad v \in M.
$$
It follows from this that $\psi_1(v)=\psi_2(\alpha v)=0$ for $v \in
M$, and hence
 $\widehat\mu_1(y)=\widehat\mu_2(\alpha y)=1, \ y \in M$.
Put $B=M\cap \alpha(M)$. Then $B$ is the required subgroup. Lemma 2
is proved.  $\Box$

\bigskip

{\bf Lemma 3}  (\cite[\S 2]{Fe6})  {\it Let $X$ be a second
countable locally compact Abelian group, and  $\mu \in {M^1}(X)$.
Let $E=\{y\in Y:\widehat \mu(y)=1\}$. Then  $E$ is a closed subgroup
of $Y$, the characteristic function $\widehat \mu(y)$ is
$E$-invariant, i.e. $\widehat\mu(y+h) = \widehat\mu(y)$ for all $y
\in Y$,  $h \in E$, and $\sigma(\mu)\subset A(X, E)$}.

\bigskip

An  Abelian group $G$ is called $p$-prime if the order of every
element of $G$ is a power of $p$. Denote by ${\cal P}$ the set of
prime numbers. The following result follows from the proof of
Theorem 1 in \cite{FeGra1} (see also \cite[\S 13]{Fe6}).

\bigskip

{\bf Lemma 4}  {\it Let $X$ be a  group of the form
$$
\mathop{\mbox{\rm\bf P}}\limits_{p \in {\cal P}}(\Delta_p^{n_p}
\times G_p),
$$
where $n_p$ is  a nonnegative integer, and  $G_p$ is a finite
$p$-prime group, may be $G_p=\{0\}$. Let $\xi_1$ and $\xi_2$ be
independent random variables with  values in $X$ and
  distributions $\mu_1$ and $\mu_2$. If the linear forms
 $L_1
= \xi_1 + \xi_2$ and $L_2=\xi_1 + \alpha \xi_2$, where $\alpha \in
{\rm Aut}(X),$ are independent, then  $\mu_j=m_K*E_{x_j}$, where $K$
is a compact subgroup of
 $X$, and $x_j\in X, \ j=1, 2$.
}

\bigskip

{\bf Lemma 5} (\cite[\S 13]{Fe6})  {\it Let $X$ be a second
countable locally compact Abelian group, $\xi_1$ and $\xi_2$ be
independent identically distributed random variables with values in
$X$ and distribution
 $m_K$, where $K$ is a compact subgroup of  $X$. Let
 $\alpha \in {\rm Aut}(X).$
Then the following statements are equivalent:

$(i)$ the linear forms  $L_1=\xi_1+\xi_2$ and $L_2=\xi_1 +
\alpha\xi_2$ are independent;

$(ii)$ $(I-\alpha)(K) \supset K.$}

\bigskip

{\it Proof of Theorem 1}  Let $\xi_1$ and $\xi_2$ be independent
random variables with values in $X$ and distributions $\mu_1$ and
$\mu_2$. Assume that the linear forms  $L_1 = \xi_1 + \xi_2$ and
$L_2=\xi_1 + \alpha \xi_2$  are independent. By Lemma 1 the
characteristic functions of the distributions $\mu_j$  satisfy
equation (\ref{3}). It is obvious that the characteristic functions
of the distributions $\bar\mu_j$ also satisfy equation (\ref{3}).
This implies that the characteristic functions of the distributions
 $\nu_j=\mu_j*\bar\mu_j$ satisfy equation (\ref{3}) as well.  We have
$\widehat\nu_j(y)=|\widehat\mu_j(y)|^2 \ge 0,$ $j = 1, 2$. Hence,
when we prove Statements  { $1(i)$} and { $1(ii)$} we may assume
without loss of generality that
 $\mu_j(y)
\ge 0,$ $j = 1, 2,$ because $\mu_j$ and $\nu_j$ are either
degenerate distributions or idempotent distributions simultaneously.
 Moreover, if it is necessary we can
consider new independent random variables  $\xi'_1=\xi_1$ and
$\xi'_2=\alpha\xi_2$, and hence, we can assume that $k \ge 0$. Note
also that the only nonzero proper closed subgroups of $\Omega_p$, are the
subgroups $\Lambda_k=p^k\Delta_p$, $k=0, \pm 1, \dots$
\cite[(10.16)]{HeRo}).

\bigskip

{\it Statement $1(i)$} We can assume that $\alpha\ne I$.
In the opposite case, obviously, $\mu_1$ and $\mu_2$ are degenerate
distributions. Since by the condition $k=0$, we have $\alpha=c,$ $c
\in \Delta_p^0$, and hence, the restriction of the automorphism
$\alpha \in {\rm Aut}(X)$ to any subgroup $p^m\Delta_p$ is a
topological automorphism of $p^m\Delta_p$. By Lemma 2 there exists a
subgroup
 $B=p^l\Delta_p$ such that $\widehat\mu_j(y)=1,$ $
j=1, 2,$ for $y\in B$.  It follows from Lemma 3 that
$\sigma(\mu_j)\subset A(X, B).$ Put $G=A(X, B)$. It is easy to see
that $G= p^{-l+1}\Delta_p$. We have $G \cong \Delta_p$, and the
restriction of $\alpha$ to the subgroup $G$ is a topological
automorphism of $G$. Thus, we get that the independent random
variables $\xi_1$ and $\xi_1$ take values in a group
$G\cong\Delta_p$, they have
 distributions $\mu_1$ and $\mu_1$, and the linear forms $L_1 = \xi_1 + \xi_2$ and
$L_2=\xi_1 + \alpha \xi_2$, where $\alpha \in {\rm Aut}(G)$, are
independent. Applying  Lemma 4, and taking into account that
$\mu_j(y) \ge 0,$ $j = 1, 2,$ we obtain that $\mu_1=\mu_2=m_K$,
where $K$ is a compact subgroup of
 $G$.  Thus, we proved the first part of Statement
 $1(i)$. On the other hand, we have independent identically distributed random
variables $\xi_1$ and $\xi_2$ with values in  $X$ and  distribution
 $m_K$  such that the linear forms $L_1=\xi_1+\xi_2$ and $L_2=\xi_1 +
\alpha\xi_2$ are independent. Hence, by Lemma 5 $(I-\alpha)(K)
\supset K.$ Suppose that $c=(0, 0, \dots, 0, 1, c_1, \dots)$,
 and  $K\ne\{0\}$. It is obvious that in this case
$(I-\alpha)(K)$ is a proper subgroup of $K$. The obtained
contradiction shows that $K=\{0\}$, i.e. $\mu_1$ and $\mu_2$ are
degenerate distributions. We also proved the second part of Statement $1(i)$.

In particular, it follows from this reasoning that in the case, when
$X=\Omega_2$, $\mu_1$ and $\mu_2$ are degenerate distributions,
because if $c \in \Delta_2^0$, then $c_0=1$.

\bigskip

{\it Statement $1(ii)$} Put $f(y)=\widehat\mu_1(y),$
$g(y)=\widehat\mu_2(y)$. Taking into account that
$\alpha=\widetilde\alpha$, we rewrite equation  (\ref{3}) in the form
\begin{equation}\label{5}
f(u+v)g(u+\alpha v)=f(u)g(u)f(v)g(\alpha v), \quad u, v \in Y.
\end{equation}
Put
\begin{equation}\label{8}
E = \{ y \in Y: f(y) = g(y) = 1 \}.
\end{equation}
Obviously, we can assume that $\mu_j$ are nondegenerate
distributions, and hence $E \ne \Omega_p$. By Lemma 2 $E \ne
\{0\},$ and by Lemma 3, $E$ is a closed subgroup of
 $\Omega_p.$ Thus, $E$, as a nonzero proper closed
subgroup of
 $\Omega_p,$ is of the form
  $E=p^l\Delta_p$. Since $k\ge 1$, we have $\alpha(E) \subset E$
and hence, $\alpha$   induces a continuous endomorphism
$\widehat\alpha$ on the factor-group $L=Y/E$. Taking into account
that by Lemma 3
$$
f(y + h) = f(y), \quad g(y + h) = g(y),
$$
for all $y \in Y,$ $h \in E$, we can consider the functions
$\widehat f(y)$ and $\widehat g(y)$ induced on $L$ by the functions
$f(y)$ and $g(y)$. It follows from (\ref{8}) that
\begin{equation}\label{9}
\{ y \in L: \widehat f(y) = \widehat g(y) = 1 \} = \{0\}.
\end{equation}
Passing from equation (\ref{5}) on the group $Y$ to the induced
equation on the factor-group $L = Y/E$, we obtain
\begin{equation}\label{5a}
\widehat f(u+v)\widehat g(u+\widehat\alpha v)=\widehat
f(u)\widehat g(u)\widehat f(v)\widehat g(\widehat\alpha v), \quad
u, v \in L.
\end{equation}
It is easy to see  that  $L \cong{\mathbb Z}(p^\infty)$ and
$\widehat\beta=(I-\widehat\alpha)\in {\rm Aut}(L)$. Putting in
(\ref{5a}) first $u = -\widehat\alpha  y$, $v = y$, and then
 $u = y,$ $v
= - y$, and taking into account that $\widehat f(-y)=\widehat f(y)$
and  $\widehat g(-y)=\widehat g(y)$, we get
\begin{equation}\label{6}
\widehat f((I-\widehat\alpha)y)=\widehat f(\widehat\alpha y)\widehat
g^2(\widehat\alpha y)\widehat f(y), \quad y \in L,
\end{equation}
\begin{equation}\label{7}
\widehat g((I-\widehat\alpha)y)=\widehat f^2(y)\widehat g(y)\widehat
g(\alpha y), \quad y \in L.
\end{equation}

Obviously,   equation (\ref{6})  implies that
\begin{equation}\label{10}
\widehat f(\widehat\beta y) \le \widehat f(y), \quad y \in L.
\end{equation}
We note now that any element of the group $L$ belongs to some
subgroup $H$, $H \cong{\mathbb Z}(p^m)$, moreover,
$\widehat\beta(H)=H$. Since $H$ is a finite subgroup,
$\widehat\beta ^ny = y$ for any $y\in H$, where
 $n$ depends generally on $y$. Then (\ref{10})  implies that
$$
\widehat f(y) = \widehat f(\widehat\beta^ny) \le \dots \le \widehat
f(\widehat\beta y) \le \widehat f(y), \quad y \in L.
$$ Thus, on each orbit
$\{y, \ \widehat\beta y, \dots, \widehat\beta^{n-1}y \}$ the
function $\widehat f(y)$ takes a constant value. The similar
statement for the function $\widehat g(y)$ follows from the equation
induced by equation (\ref{7}).

Assume that $\widehat f(y_0) \ne 0$  at a point  $y_0 \in L,$ $y_0
\ne 0$.
 Then $\widehat f(\widehat\beta y_0)
= \widehat f(y_0) \ne 0$, and    equation (\ref{6}) implies that
\begin{equation}\label{11}
\widehat f(\widehat\alpha y_0) = \widehat g(\widehat\alpha y_0) = 1.
\end{equation}

It follows from (\ref{9}) and (\ref{11}) that $\widehat\alpha y_0 =
0$. By the condition $\alpha=p c,$ where $c \in \Delta_p^0$. This
implies that $\widehat\alpha=p \widehat c$, where $\widehat c$ is an
automorphism of the group
 $L$, induced by the automorphism  $c$.
Hence, $y_0$ is an element of order $p$. Reasoning similarly we get
from   equation
 (\ref{7}) that if  $\widehat g(y_1) \ne 0,$ $y_1 \in L,$ $y_1
\ne 0$, then $\widehat f(y_1) = \widehat g(\widehat\alpha y_1) = 1$.

Let $w$ be an arbitrary element of $L$. Denote by $\langle w
\rangle$ the subgroup of $L$ generated by $w$. It follows from
$\widehat f(y_1)=1$ that  $\widehat f(y)=1$ for all $y \in \langle
y_1 \rangle$. Since
 $L \cong{\mathbb
Z}(p^\infty)$ and  $\langle y_1 \rangle$ is a subgroup of $L$, we
have $\langle y_1 \rangle\cong{\mathbb Z}(p^m)$ for some  $m$, and
hence
 $\widehat\alpha (\langle y_1 \rangle)\subset \langle y_1 \rangle$. Moreover,
$\widehat f(\widehat\alpha y_1)=1$. Thus the equalities
\begin{equation}\label{12}
\widehat f(\widehat\alpha y_1) = \widehat g(\widehat\alpha y_1) = 1
\end{equation}
hold true. It follows from (\ref{9}) and (\ref{12}) that
$\widehat\alpha y_1 = 0$, and hence
 $y_1$ is also an element of order $p$.
Since $L \cong{\mathbb Z}(p^\infty)$, the group  $L$ contains the
only subgroup  $A$ topologically isomorphic to ${\mathbb Z}(p)$. So,
we proved that the functions
 $\widehat f(y)$ and $\widehat g(y)$ vanish for $y \notin A$.

Consider the restriction of equation (\ref{5a}) to the subgroup $A$.
Taking into account that $\widehat\alpha y = 0$ for all $y \in A$,
we obtain
\begin{equation}\label{13}
\widehat f(u+v)\widehat g(u)=\widehat f(u)\widehat g(u)\widehat
f(v), \quad u, v \in A.
\end{equation}

If $\widehat g(u_0) \ne 0$ at a point $u_0 \in A,$ $u_0 \ne 0$, then
we conclude from (\ref{13}) that
$$
\widehat f(u_0+v)=\widehat f(u_0)\widehat f(v), \quad v \in A.
$$
Putting here $v = (p-1)u_0$, we get $\widehat f(u_0) = 1$. Since $p$
is a prime number, we have
 $A =\langle u_0 \rangle$, and hence, $\widehat f(y) = 1$ for $y \in A$. If $\widehat
g(y) = 0$ for any $y \in A,$ $y \ne 0$, then, obviously, $\widehat
f(y)$ may be an arbitrary positive definite function on $A$. Thus,
we proved that either
\begin{equation}\label{14}
\widehat f(y) = \cases{ 1, \quad  \quad y \in A, \cr\cr 0, \quad
\quad y \notin A, \cr}
\end{equation}
or
\begin{equation}\label{15}
\widehat g(y) = \cases{ 1, \quad  \quad y =0, \cr\cr 0, \quad  \quad
y \ne 0. \cr}
\end{equation}
Return  from the induced functions $\widehat f(y)$ and $\widehat
g(y)$ on $L$ to the functions $f(y)$ and $g(y)$ on $Y$. Taking
into account (\ref{2}) and the fact that a distribution is
uniquely defined by its characteristic   function, we obtain from
(\ref{14}) and (\ref{15}) that either $\mu_1 \in I(X)$, or $\mu_2
\in I(X)$. Statement    $1(ii)$ is proved.

\bigskip

{\it Statement $2$} It is easy to see that without loss of generality we
can assume that $k\ge 2$. Consider on the group $\Omega_p$ the
distributions $$\mu_1=a m_{\Lambda_{1}}+(1-a)m_{\Lambda_{-k+2}},
\quad \mu_2=a m_{\Lambda_{-k+2}}+(1-a)m_{\Lambda_{-k+1}},$$ where
$0<a<1$. As has been noted earlier, $A(Y,
\Lambda_m)=\Lambda_{-m+1}$. Therefore (\ref{2})  implies that the
characteristic functions $f(y)=\widehat\mu_1(y)$ and
$g(y)=\widehat\mu_2(y)$ are of the form
$$
f(y) = \cases{ 1, \quad\quad  \quad y \in p^{k-1}\Delta_p, \cr\cr
a,\quad\quad y \in \Delta_p\backslash p^{k-1}\Delta_p, \cr\cr 0,
\quad\quad\quad\quad  \quad y \notin\Delta_p, \cr}
\quad\quad\quad\quad g(y) = \cases{ 1, \quad\quad\quad\quad  \quad y
\in p^{k}\Delta_p, \cr\cr a,\quad\quad y \in
p^{k-1}\Delta_p\backslash p^{k}\Delta_p, \cr\cr 0,
\quad\quad\quad\quad  \quad y \notin p^{k-1}\Delta_p. \cr}
$$
Let $\xi_1$ and $\xi_2$ be independent random variables with values
in the group
 $\Omega_p$ and distributions $\mu_1$ and $\mu_2$. It is obvious that $\mu_1, \mu_2 \notin I(X)$.
 We will check that
the characteristic functions  $f(y)$ and $g(y)$ satisfy equation
(\ref{5}). Then, by Lemma 1, the linear forms $L_1=\xi_1 +  \xi_2$
and $L_2= \xi_1 + \alpha\xi_2$ are independent, and Statement 2
will be proved.

\bigskip

Consider 3 cases: {\bf 1}. $u, v \in \Delta_p$; {\bf 2}. $u \notin
\Delta_p, \  v \in \Delta_p$; and {\bf 3}. $v \notin \Delta_p$.

\bigskip

{\bf 1}. $u, v \in \Delta_p$. Note that since $k\ge 2$, we have
$p^{k-1}\Delta_p\subset \Delta_p$. Consider  3 subcases.

\bigskip

${\bf  1a}$. $u \in p^{k-1}\Delta_p,$ $v \in \Delta_p$. Since
 $u \in p^{k-1}\Delta_p$, we have $f(u)=1$, and hence
$f(u+v)=f(v)$. Since $\alpha v \in p^{k}\Delta_p$, we have
$g(\alpha v)=1$, and hence $g(u + \alpha v)=g(u)$. Equation
(\ref{5}) takes the form $f(v)g(u)=f(v)g(u)$, and it is obviously
true.

\bigskip

${\bf 1b}$. $u \in \Delta_p\backslash p^{k-1}\Delta_p, \ v \in
p^{k-1}\Delta_p.$ Since $v \in p^{k-1}\Delta_p,$ we have $\alpha v
\in p^{2k-1}\Delta_p \subset p^{k}\Delta_p.$ This implies that
$g(\alpha v)=1$, and hence, $g(u + \alpha v)=g(u)$. Since $v \in
p^{k-1}\Delta_p,$ we have
 $f(v)=1$, and hence
$f(u+v)=f(u).$ Equation (\ref{5}) takes the form
$f(u)g(u)=f(u)g(u)$, and it is obviously true.

\bigskip

${\bf 1c.}$ $u \in \Delta_p\backslash p^{k-1}\Delta_p, \ v \in
\Delta_p\backslash p^{k-1}\Delta_p.$ Since $v \in \Delta_p$, we have
$\alpha v \in p^{k}\Delta_p.$ This implies that   $g(\alpha v)=1$,
and hence $g(u + \alpha v)=g(u)$. Since  $u \notin p^{k-1}\Delta_p,$
we have
 $g(u + \alpha v)=g(u)=0$. Thus, both sides of equation (\ref{5}) vanish.

\bigskip

{\bf 2.} $u \notin \Delta_p, \  v \in \Delta_p$. This implies that
$u + v \notin \Delta_p$, and hence $f(u)=0$ and $f(u+v)=0$. Thus
both sides of equation (\ref{5}) vanish.

\bigskip

{\bf 3}. $v \notin \Delta_p$. This implies that $f(v)=0$ and hence,
the right-hand side of equation (\ref{5}) vanishes. If the left-hand
side of equation (\ref{5}) does not vanish, then the following
inclusions
\begin{equation}\label{15a}
\cases{ u+v \in \Delta_p, \cr\cr u+\alpha v \in p^{k-1}\Delta_p\cr}
\end{equation}
hold true. On the one hand, since $k\ge 2$, it follows from
(\ref{15a}) that $(I-\alpha)v \in \Delta_p$. On the other hand,
since $k\ge 2$, we have $(I-\alpha) \in {\rm Aut}(\Delta_p)$. Hence
$v \in \Delta_p$. The obtained contradiction shows that the
left-hand side of equation (\ref{5})   vanishes as well.

We showed that the characteristic functions  $f(y)$ and $g(y)$
satisfy equation (\ref{5}). Thus, we proved Statement 2 and hence,
Theorem 1 is completely proved. $\Box$

\bigskip

{\bf Remark 1}  As follows from the proof of Statement
 $1(i)$ if $k=0$, then
 $\mu_j=m_K*E_{x_j}$, where $K$ is a compact subgroup of
 $\Omega_p$, $x_j\in
 \Omega_p$, $j=1, 2$.

\bigskip

As a corollary from Theorem 1 and Remark 1 we derive  the
Kac-Bernstein theorem for the group $\Omega_p$ (see \cite[\S
7]{Fe6}).

\bigskip

 {\bf Corollary
1}  {\it Let   $\xi_1$ and $\xi_2$ be independent random variables
with values in $\Omega_p$ and distributions $\mu_1$ and $\mu_2$.
Assume that the linear forms $L_1 = \xi_1 + \xi_2$ and $L_2=\xi_1
- \xi_2$ are independent. If $p=2$, then $\mu_1$ and $\mu_2$ are
degenerate distributions. If $p>2$, then $\mu_j=m_K*E_{x_j}$,
where $K$ is a compact subgroup of
 $\Omega_p$, $x_j\in
 \Omega_p$, $j=1, 2$.}

\bigskip

{\bf Remark 2}  Let $\xi_1$ and $\xi_2$ be independent random
variables with values in the group $\Omega_p$ and  distributions
$\mu_1$ and $\mu_2$. Assume that the linear forms $L_1=\alpha_1
\xi_1 + \alpha_2 \xi_2$ and $L_2=\beta_1 \xi_1 + \beta_2 \xi_2,$
where $\alpha_j, \ \beta_j \in {\rm Aut}(\Omega_p)$, are
independent. We can consider new independent random variables
 $\xi'_1 =
\alpha_1\xi_1$ and $\xi'_2 = \alpha_2\xi_2$ and reduce the problem
of describing possible distributions $\mu_1$ and $\mu_2$ to the
case, when $L_1= \xi_1 + \xi_2, \ L_2=\delta_1 \xi_1 + \delta_2
\xi_2,$ where $\delta_1, \ \delta_2 \in {\rm Aut}(\Omega_p)$.
Since $L_1$ and $L_2$ are independent if and only if  $L_1$ and
$L'_2=\delta_1^{-1} L_2$ are independent, the problem of
describing possible distributions $\mu_1$ and $\mu_2$ is reduced
to the case when  $L_1= \xi_1 + \xi_2$ and $L_2= \xi_1 + \alpha
\xi_2$, where $\alpha \in {\rm Aut}(\Omega_p)$, i.e. it is reduced
to Theorem 1.

\bigskip

{\bf Remark 3} Consider the group  $\Omega_p,$ where $p>2$. Let
$\xi_1$ and $\xi_2$
 be independent identically distributed
random variables with  values in $\Omega_p$ and
  distribution $m_{\Delta_p}$.  Let
$\alpha=(0, 0, \dots, 0, x_0, x_1, \dots)\in {\rm Aut}(\Omega_p)$,
where  $x_0\ne 1$. It is easy to verify that the characteristic
functions $\widehat\mu_1(y)=\widehat\mu_2(y)=\widehat
m_{\Delta_p}(y)$ satisfy equation (\ref{3}). This implies by Lemma 1
that the linear forms
   $L_1 = \xi_1
+ \xi_2$ and $L_2=\xi_1 + \alpha \xi_2$ are independent. Thus, for
the group $\Omega_p,$ where $p>2$, Statement $1(i)$ can not be
strengthened to the statement that  both $\mu_1$ and $\mu_2$ are
degenerate distributions.

\bigskip

{\bf Remark 4} Statement $1(ii)$ can not be strengthened to the
statement that both $\mu_1$ and $\mu_2$ are idempotent
distributions. Namely, if $k=1$, then there exist independent
random variables $\xi_1$ and $\xi_2$  with  values in the group
$X=\Omega_p$ and distributions $\mu_1$ and $\mu_2$ such that the
linear forms
   $L_1 = \xi_1
+ \xi_2$ and $L_2=\xi_1 + \alpha \xi_2$ are independent, but one of
the distributions
 $ \mu_j \notin I(X)$. We  get the corresponding example if we put
$\mu_1= m_{\Lambda_{1}}$ and $\mu_2=a
m_{\Lambda_{1}}+(1-a)m_{\Lambda_{0}}$, where $0<a<1$. The proof is
similar to the reasoning given in the proof of Statement 2.

\end{document}